\documentclass[12pt]{article}
\usepackage{amsmath,amsfonts,amssymb}
\usepackage{secdot}
\usepackage{amsthm}
\usepackage{color}
\usepackage[russian]{babel}
\newtheorem{theorem}{\hspace*{\parindent}Theorem}

\newcounter{theremark}

\newcommand{\rr}{\mathbb R^d}
\def\a{\mathbf{a}}
\def\x{\mathbf{x}}
\def\y{\mathbf{y}}
\title{Об оптимальной дискретной энергии  Неймана \\ в шаре и круговом кольце}
\author{А. С. Афанасьева-Григорьева, Е. Г. Прилепкина}
\date{}

\begin{document}
\maketitle

\begin{abstract}
В работе доказаны некоторые точные оценки дискретной энергии Неймана шара и кругового кольца в евклидовом пространстве для точек, расположенных на окружностях. Доказательства основаны на диссимметризации и анализе асимптотического поведения интеграла Дирихле потенциальной функции. 
\end{abstract}

\bigskip

\emph{Ключевые слова:} дискретная энергия, функция Грина, функция Неймана,
диссимметризация

\bigskip


\bigskip

\section{Введение и формулировки результатов}

 В данной работе  $\mathbb{R}^{d}$ будет означать  $d$-мерное евклидово пространство точек  $\x$ вида $(x_1,\dots,\,x_d)$ с обычной длиной и расстоянием, $d\geq 2$. В случае $d=2$ мы считаем, что $\mathbb R^2$ является комплексной плоскостью. 
 Решение классической проблемы Неймана в  ограниченной области $D\subset \mathbb{R}^{d}$ для уравнения Пуассона требует построения  функции Неймана (иногда ее называют функцией Грина для проблемы Неймана или функцией Грина второго рода).  Классическая  функция Неймана  определяется \cite{Henr}, \cite{Sad}  как функция $\x\in D$ в области $D\setminus\{\y\}$,  имеющая представление 
\begin{equation}\label{def1}
N(\x,\y,D)=\frac{\mu_d(|\x-\y|)+v(\x,\y,D)}{w_d} 
\end{equation}
и удовлетворяющая условиям $$\frac{\partial N(\x,\y,D)}{\partial n_\x}=-\frac{1}{s_{d-1}(\partial D)},$$
$$\int_{\partial D} N(\x,\y,D)d\sigma_\x=0.$$
Здесь  $\mu_d(\cdot)$---  фундаментальное решение уравнения Лапласа,   ($\mu_2(\rho)=-\log \rho,$ 
$\mu_d(\rho)={\rho^{2-d}}/{(d-2)}$ при $\ d\geq 3$),  $w_d=2\pi^{d/2}/\Gamma(d/2)$ --- площадь единичной гиперсферы, $v(\x,\y,D)$---  некоторая гармоническая  в области $D$ функция, $s_{d-1}$--- мера Лебега  и дифференцирование берется по внешней  нормали. 

Существует много исследований, связанных с экстремальными задачами для различных видов энергий дискретного заряд (см., например, работы \cite{Brauchart}, \cite{Borodachov}, \cite{Brauchart2} и ссылки в них).   В \cite{Dubinin1} получены две оценки дискретной энергии  функции Грина кругового кольца на плоскости в случае  точек, расположенных на некоторой окружности. Эти результаты были распространены в  евклидово пространство в  \cite{DubPrilGreen}. Целью настоящей работы является получение результатов подобного сорта для функции Неймана.  

 Напомним определение дискретной энергии Грина  \cite{Lan}. Пусть  $\Delta=\{\delta_k\}_{k=1}^n$ произвольный дискретный заряд 
(множество вещественных чисел), принимающий значение  $\delta_k$ в точке  $\mathbf x_k$, $k=1,\ldots,n$ области $D$ . \emph{Энергией Грина} этого заряда относительно области $D$ называется  величина
$$
E(X,\Delta,D)=\sum\limits_{k=1}^n\sum \limits_{{l=1} \atop {l\not=
k}}^n \delta_k\delta_l g_{\small D}(\mathbf x_k,\mathbf x_l),
$$
где $g_{\small D}(\mathbf x_k,\mathbf x_l)$ функция Грина области D. Аналогичным образом определим  \emph{энергию Неймана}
$$
En(X,\Delta,D)=\sum\limits_{k=1}^n\sum \limits_{{l=1} \atop {l\not=
k}}^n \delta_k\delta_l N(\mathbf x_k,\mathbf x_l,D).
$$
 
  Всюду далее  область $D$ является либо шаром вида $\{|\x|<\tau\}$, либо  концентрическим круговым кольцом вида $\{\tau_1<|x|<\tau_2\}$.  
 Примем следующие обозначения: $B(\a,r)$ --- открытый  шар с центром в точке $\a$ радиуса $r$,  $J$--   $(d-2)$--мерная плоскость $\{\mathbf
x\in\rr:\mathbf x=(0,0,x_3,\ldots, x_d)\}$. Нам понадобятся цилиндрические координаты $(r,\theta, \mathbf{x}')$ точки $\mathbf{x}=(x_1,\dots,\,x_d)$ в  $\rr$,  связанные с декартовыми координатами соотношениями  $x_1=r \cos \theta,$ $x_2=r \sin \theta,$
$\mathbf{x}'\in J$. Записи типа $\{\theta=\varphi\}$ означают множество точек $\mathbb{R}^{d}$, имеющих   полярные координаты $(r,\varphi,x'),$  $r\geq 0$, $\x'\in J$, $\varphi$ фиксировано.

Пусть 
$\Omega=\{S\}$ означает множество, состоящее из конечного числа различных  окружностей    $S$ вида  $S=\{(r_0,\theta,\mathbf
x'_0):0\leq\theta\leq 2\pi\}$, лежащих в области $D$ (здесь 
$r_0>0$ и  $\mathbf x'_0\in J$ предполагается фиксированным). Для произвольных вещественных чисел $\theta_j,$ $j=0,\ldots,$ \makebox{$m-1,$}
\begin{equation*}
0\leq \theta_0<\theta_1<\ldots<\theta_{m-1}<2\pi,
\end{equation*}
обозначим $X=\{\mathbf x_k\}_{k=1}^n$ множество точек пересечения окружностей из $\Omega$
с полуплоскостями 
\begin{equation*}
L_j=\{(r,\theta,\mathbf x'):\theta=\theta_j\}, \ j=0,\ldots,m-1.
\end{equation*}
Обозначим также  $X^*=\{\mathbf x^*_k\}_{k=1}^n$ -- множество точек пересечения окружностей из  $\Omega$ 
с симметричными полуплоскостями 
\begin{equation*}
\ L_j^*=\{(r,\theta,\mathbf x'): \theta=2\pi j/m\}, \
j=0,\ldots,m-1.
\end{equation*}

Следующие теоремы показывают, что в зависимости от условий на заряд $\Delta$ симметричная конфигурация дает как максимум, так и минимум энергии Неймана $En(X,\Delta,D)$. 

\begin{theorem}
Пусть $D$ шар или круговое кольцо, $\Omega$, $X$ и 
$X^*$ определены выше,  заряд  $\Delta=\{\delta_k\}_{k=1}^n$  принимает одинаковые значения $\delta_k=\delta_l$ в точках  $\mathbf x_k\in X$ и 
$\mathbf x_l\in X$, расположенных на одной и той же окружности 
из  $\Omega$ и $$\sum_{k=1}^n\delta_k=0. $$ Кроме того, пусть точки  $\mathbf x_k\in X$
и  $\mathbf x^*_k\in X^*$ лежат на одной и той же окружности из $\Omega$,
$k=1,\ldots,n$.  Тогда 
$$
En(X,\Delta,D)\geq En(X^*,\Delta,D).
$$
\end{theorem}

\begin{theorem}
Пусть $D$ шар или круговое кольцо, $\Omega$, $X$,  
$X^*$, $\Delta$ определены выше,   $m$ --- четное число и $\delta_k=-\delta_l$  в точках 
$\mathbf x_k\in X$ и $\mathbf x_l\in X,$ лежащих на одной и той же окружности из $\Omega$ и  на соседних полуплоскостях из совокупности $\{L_j\}_{j=0}^{m-1}$. Тогда 
$$
En(X,\Delta,D)\leq En(X^*,\Delta,D),
$$
где точки   $X^*$ пронумерованы следующим образом:
если $\mathbf x_k^*\in X^*$ лежит на пересечении окружности  $S$
из   $\Omega$ с полуплоскостью  $L_j^*$, тогда соответствующая точка $\mathbf x_k\in X$ должна лежать на пересечении  $S$ и полуплоскости  $L_j$, $k=1,\ldots,n$,  $0\leq j\leq
m-1$.
\end{theorem}

 Заметим, что полученные в работе теоремы справедливы и в случае, когда $D$ означает область вращения (область  $D\subset \rr$ называется 
\emph{областью вращения} относительно оси  $J$, если для любой точки  $(r,\theta,\mathbf{x}')\in B$ и любого  $\varphi$ точка  $(r,\varphi,\mathbf{x}')$ принадлежит  $D$). 

При дополнительном   условии

\begin{equation}\label{condit}
\sum_{k=1}^{n}\delta_{k}=0
\end{equation}
определим функцию 
\[
u(\x)=u(\x;X,D,\Delta)=\sum_{k=1}^{n}\delta_{k}N(\x,\x_{k},D),
\]
которую назовем потенциальной функцией Неймана конфигурации $X,$ $\Delta$, $D$.  Непосредственно из определения  вытекает разложение потенциальной функции в окрестности точки $\x_{k},$ $k=1,\dots,n,$
\begin{equation}\label{eq:AK}
u(\x)=\delta_{k}\frac{\mu_d(|\x-\x_k|)}{w_d}+a_{k}+o(1), \, \x\to\x_k,
\end{equation}
где 
\[
a_{k}=\delta_{k}\frac{v(\x_k,\x_k,D)}{w_d}+ \sum_{\substack{l=1\\
l\ne k
}
}^{n}\delta_{l}N(\x_{l},\x_{k},D).
\]

Сумма 
\begin{equation}\label{quadratic}\sum_{k=1}^n\delta_k a_k=\sum_{k=1}^n\sum_{l=1}^n
\eta_{kl}(D)\delta_k\delta_{l}=En(X,\Delta,D)+\sum_{k=1}^{n}\frac{\delta_k^2 v(\x_k,\x_k,D)}{w_d}
\end{equation} 
представляет из себя квадратичную форму переменных $\Delta$ с коэффициентами $\eta_{kl}(D)$, зависящими от функции Неймана.
Обозначим эту квадратичную форму 
\begin{equation}\label{QNeim}
   Qn(X,\Delta,D)=\sum_{k=1}^n\sum_{l=1}^n
\eta_{kl}(D)\delta_k\delta_{l},
\end{equation}
где $\eta_{kl}(D)=N(\x_k,\x_l,D),\ k\neq l,$ $\eta_{kk}(D)=v(\x_k,\x_l,D)/w_d.$

Квадратичные формы такого сорта, а также формы с коэффициентами, зависящими от функции Грина либо Робена, играют важную роль в геометрической теории функций.  Различные неравенства для таких форм и их применения встречаются в работах Аленицина, Нехари, Дюрена, Шиффера, Дубинина и других математиков (см.  \cite{N}, \cite{Duren}, \cite{Dub7}, \cite{PrilPomi}).  Мы доказываем, что 
$$Qn(X,\Delta,D)\geq Qn(X^*,\Delta,D)$$ в условиях    Теоремы 1 и  
$$Qn(X,\Delta,D)\leq Qn(X^*,\Delta,D)$$ 
 в условиях Теоремы 2. Для вычисления коэффициентов квадратичной формы $Q_n$ при дополнительном условии \eqref{condit} взамен классической можно использовать обобщенную функцию Неймана \cite{KarpPrilepkina}.  На плоскости известен  явный вид формы  $Q_n$ круга  и  кольца.  Функция Неймана единичного круга $U$ \cite{Sad}
  \[
N\left(z,z_0, U\right)=-\frac{\log|z-z_0|1-z\overline{ z_0}|}{2\pi},\]   поэтому 
$$\eta_{kl}(U)=-\frac{\log|z_k-z_l|1-z_k\overline{ z_l}|}{2\pi}, \ k\neq l,$$
$$\eta_{kk}(U)=-\frac{\log(1-|z_k|^2)}{2\pi}.$$
В \cite{KarpPrilepkina} приведены коэффициенты $\eta_{kl}(K)$ квадратичной формы плоского кольца $K=\{\mu<|z|<1\}.$ А именно, 
$$
\eta_{kl}(K)=
\begin{cases}
-\frac{1}{2\pi}\log{|\theta_1(i\log(z_{k}\overline{z_{l}})/2;\mu)\theta_1(i\log(z_{k}/z_{l})/2;\mu)|},& k \not= l,\\
\frac{1}{2\pi}\log\frac{4|z_k|^2|\sin{(i\log{|z_k|}})|}{(1-|z_k|^2)|\theta_1(i\log|z_k|;\mu)\theta_1'(0;\mu)|},& k=l,
\end{cases}
$$
где
$$
\theta_1(z;\mu)=-i\sum_{n=-\infty}^{\infty}(-1)^{n}\mu^{(n+1/2)^2}e^{i(2n+1)z}.$$

В пространстве размерности $d\geq 3$ мы не нашли в литературе аналитического выражения функции Неймана кругового  кольца. Для единичного шара  $U=B(0,1)$  функция Неймана найдена в работе \cite{Sad} и имеет вид   
$$
N(\x,\y,U)=\frac{1}{\omega_d}\left(\mu_d(|\x-\y|)+\mu_d\left(\left|x|\y|-\frac{\y}{|\y|}\right|\right)+\epsilon_1(\x,\y)\right)+Const,
$$
где  $\epsilon_1(\x,\y)$ задается формулами 

$$
\epsilon_1(\x,\y)=\log{\frac{2}{\left|1-(\x,\y)+\left|\x|\y|-\frac{\y}{|\y|}\right|\right|}}, d=3;
$$

$$
\epsilon_1(\x,\y)=\frac{(\x,\y)}{\sqrt{|\x|^{2}|\y|^{2}-(\x,\y)^2}}\arctan\frac{\sqrt{|\x|^{2}|\y|^{2}-(\x,\y)^2}}{1-(\x,\y)}-\log{\left|\x|\y|-\frac{\y}{|\y|}\right|}, d=4;
$$

\begin{multline*}
\epsilon_1(\x,\y)=\log{\frac{2}{\left|1-(\x,\y)+\left|\x|\y|-\frac{\y}{|\y|}\right|\right|}}+\sum\limits_{k=1}^{p-1}\frac{1}{(2k-1)}\left(\left|\x|\y|-\frac{\y}{|\y|}\right|^{1-2k}-1\right)\\
+\sum\limits_{k=1}^{p-1}\sum\limits_{i=0}^{p-k-1}\frac{2^i(k+i-1)!(2k-3)!!}{(k-1)!(2k+2i-1)!!}\frac{(\x,\y)|\x|^{2i}|\y|^{2i}}{\left(|\x|^{2}|\y|^{2}-(\x,\y)^{2}\right)^{i+1}}\left(\frac{|\x|^{2}|\y|^{2}-(\x,\y)}{\left|\x|\y|-\frac{\y}{|\y|}\right|^{2k-1}}+(\x,\y)\right),
\\d\geq 5, d=2p+1,p\geq 2;
\end{multline*}

\begin{multline*}
\epsilon_1(\x,\y)=-\log{\left|\x|\y|-\frac{\y}{|\y|}\right|}+\sum\limits_{k=1}^{p-1}\frac{1}{2k}\left(\left|\x|\y|-\frac{\y}{|\y|}\right|^{-2k}-1\right)\\
+(\x,\y)\arctan\frac{\sqrt{|\x|^{2}|\y|^{2}-(\x,\y)}^{2}}{(1-(\x,\y))}\sum\limits_{k=0}^{p-1}\frac{(2k-1)!!}{2^{k}k!}\frac{|\x|^{2k}|\y|^{2k}}{\left(|\x|^{2}|\y|^{2}-(\x,\y)^{2}\right)^{k+\frac{1}{2}}}
\\
+\sum\limits_{k=1}^{p-1}\sum\limits_{i=0}^{p-k-1}\frac{(2k+2i-1)!!(k+1)!}{2^{i+1}(2k-1)!!(k+i)!}\frac{(\x,\y)|\x|^{2i}|\y|^{2i}}{\left(|\x|^{2}|\y|^{2}-(\x,\y)^{2}\right)^{i+1}}\left(\frac{|\x|^{2}|\y|^{2}-(\x,\y)}{\left|\x|\y|-\frac{\y}{|\y|}\right|^{2k}}-(\x,\y)\right),
\\d\geq 6, d=2p+2,p\geq 2;
\end{multline*}
$0!=1,$ $(-1)!!=1.$
\section{Доказательство Теоремы 1.}

Обозначим символом $D_r$ область, полученную удалением из $D$ шаров с центром $\x_k$ радиуса $r,$ $D_r=D\setminus(\cup_{k=1}^n \overline{B(\x_k,r)}))$.  Тогда для интеграла Дирихле $I(u,D_r)=\int_{D_r}|\nabla u|^2d\x$ потенциальной функции справедлива асимптотическая  формула \cite[Лемма 2.1]{GKP}, \cite[Лемма 1]{DubPril5}
\begin{equation}\label{ass1}
I(u,D_{r})=\left(\sum_{k=1}^{n}\delta_{k}^{2}\right)\frac{\mu_d(r)}{w_d}+En(X,\Delta,D)+\sum_{k=1}^{n}\frac{\delta_k^2 v(\x_k,\x_k,D)}{w_d}+o(1),\,\,r\rightarrow0, 
\end{equation} 
или 
\begin{equation}\label{ass2}
I(u,D_{r})=\left(\sum_{k=1}^{n}\delta_{k}^{2}\right)\frac{\mu_d(r)}{w_d}+\sum_{k=1}^n\delta_k a_k+o(1),\,\,r\rightarrow0. 
\end{equation} 

Функцию  $v(\x)$ назовем допустимой для  $D$, $X$, $\Delta$,  
если $v(\x)\in\text{Lip}$ в окрестности каждой точки   $D$ за ислючением, может быть, конечного числа точек,  непрерывна в $\overline{D} \setminus \bigcup_{k=1}^{n}\{\x_{k}\}$,  и в окрестности  $\x_{k}$ справедливо разложение 
\begin{equation}\label{eq:BK}
v(\x)=\delta_{k}\frac{\mu_d(|\x-\x_k|)}{w_d}+b_{k}+o(1), \, \x\to\x_k.
\end{equation}
Для допустимой функции $v$ и потенциальной функции  $u$ мы имеем ассимптотику \cite[Лемма 2.2]{GKP}, \cite[Лемма 2]{PrilPomi}
\begin{equation}\label{eq:ass3}
I(v-u,D_{r})=I(v,D_{r})-I(u,D_{r})-2\sum_{k=1}^{n}\delta_{k}(b_{k}-a_{k})+o(1),\,\,r\rightarrow0.
\end{equation}

Пусть $u_1(\x)$ потенциальная функция Неймана набора $X$, $\Delta$, $u_2(\x)$ потенциальная функция Неймана набора $X^*$, $\Delta$, и $Dis$ означает диссимметризацию, описанную в доказательстве Теоремы 1 работы \cite{DubPrilGreen}. Построим в области $D$ функцию $v(\x)$ по правилу
$$v(\x)=u_2(Dis^{-1}(\x)).$$ В силу симметричности конфигурации $X^*$, $\Delta$, $D$ функция $u_2(\x)$ инварианта относительно любого отображения из группы симметрий $\varphi\in \Phi$, учавствующих в определении диссимметризации $Dis$. Поэтому $v(\x)$ определена однозначно и является допустимой для $X$, $\Delta$. Так как диссимметризация является, по сути, специальной перестановкой углов, то 
$$I(v,D_r)=I(u_2, D_r^*),$$
где $D_r^*=D\setminus(\cup_{k=1}^n \overline{B(\x_k^*,r)}))$. Из  \eqref{ass2}, \eqref{eq:ass3} следует 
\begin{multline}0\leq
I(v,D_{r})-I(u_1,D_{r})-2\sum_{k=1}^{n}\delta_{k}(b_{k}-a_{k})+o(1)=\\I(u_2,D_{r})-I(u_1,D_{r})-2\sum_{k=1}^{n}\delta_{k}(b_{k}-a_{k})+o(1)=\sum_{k=1}^{n}\delta_{k}(a_{k}-b_{k})+o(1),   \ \ r\to 0,
    \end{multline}
    \begin{equation}\label{prog13}
   \sum _{k=1}^{n} \delta_{k}b_{k}\leq \sum _{k=1}^{n} \delta_{k}a_{k}.
\end{equation}
Здесь $b_k$ коэффициенты асимптотического разложения потенциальной функции симметричной конфигурации, $a_k$ --- не симметричной.   С учетом  \eqref{quadratic},  получаем 
\begin{equation}\label{almost}
En(X^*,\Delta,D)+\sum_{k=1}^{n}\frac{\delta_k^2 v(\x_k^*,\x_k^*,D)}{w_d}\leq En(X,\Delta,D)+\sum_{k=1}^{n}\frac{\delta_k^2 v(\x_k,\x_k,D)}{w_d}. 
\end{equation}
Поскольку $D$ является шаром или кольцом, $v(\x,\x,D)=v(\y,\y,D)$  для любых двух точек $\x, \y,$ принадлежащих одной и той же окружности $S$ из $ \Omega$.  Следовательно,   $$\sum_{k=1}^{n}\frac{\delta_k^2 v(\x_k^*,\x_k^*,D)}{w_d}=\sum_{k=1}^{n}\frac{\delta_k^2 v(\x_k,\x_k,D)}{w_d}.$$  Таким образом, неравенство \eqref{almost} доказывает Теорему 1. 

\section{Доказательство Теоремы 2.}  
Докажем сперва вспомагательную лемму.  

{\bf{Лемма 1.}}
{\it Пусть $Y=\{\y_q\}_{q=1}^l$ совокупность точек, лежащих на полуплоскости  $\{\theta=0\},$ $\Delta_0=\{\sigma_q\}_{q=1}^l$ некоторый заряд,   $0<\alpha<\pi$,  $D(\alpha)=D\cap\{0<\theta<\alpha\}$, $\Gamma(\alpha)=\partial D(\alpha)\cap \{\theta=\alpha\}$ либо $D(\alpha)=D\cap\{-\alpha<\theta<0\}$, $\Gamma(\alpha)= \partial D(\alpha)\cap \{\theta=-\alpha\}$.
Рассмотрим функцию $h_\alpha(\x)$, гармоническую в $D(\alpha)$ за исключеним точек $Y,$ непрерывную в $\overline{D(\alpha)}\setminus Y,$  равную нулю на $\Gamma(\alpha)$,  имеющую нулевую производную на оставшейся части границы $\partial D(\alpha)\setminus Y$ и в окрестности точек $\y_q$  удовлетворяющую разложению 
\begin{equation}\label{eq:CKK}
h_\alpha(\x)=\sigma_{q}\frac{\mu_d(|\x-\y_q|)}{w_d}+c_{q}(\alpha)+o(1), \, \x\to\y_q, 
\end{equation}
Тогда  функция 
$$f(\alpha)=\sum_{q=1}^l\sigma_q c_q(\alpha)$$ вогнута на $0<\alpha<\pi$ как функция от $\alpha.$} 

{\bf Доказательство Леммы 1.} 
 Вне области $\overline{D(\alpha)}$ мы полагаем, что  функция $h_\alpha$ доопределена нулем.  В терминах работ \cite{GKP}, \cite{DubPril5} функция $h_\alpha(\x)$ называется  потенциальной функцией  набора $D(\alpha),$ $\Gamma(\alpha),$ $Y,$ $\Delta_0.$ Повторяя доказательство леммы 2.1 работы \cite{GKP} получим  разложение 
\begin{equation}\label{ass11}
I(h_\alpha,D(\alpha)_{r})=\frac{1}{2}\left(\sum_{q=1}^{l}\sigma_q^2\right)\frac{\mu_d(r)}{w_d}+\frac{1}{2}\sum_{q=1}^l\sigma_q c_q(\alpha)+o(1),\,\,r\rightarrow0. 
\end{equation} Для $0<\alpha<\beta<\pi$ построим в области $D((\alpha+\beta)/2)$  функцию $v_{(\alpha+\beta)/2}(\x)$ по правилу 
$$v_{(\alpha+\beta)/2}(\x)=\frac{h_\alpha(\x)+h_\beta(\x)-h_\beta(\x^*)}{2},$$ где $\x^*$ означает точку, симметричную $\x$ относительно полуплоскости  $\{\theta=(\alpha+\beta)/2\}$ (либо   $\{\theta=-(\alpha+\beta)/2\}$). Функция $v_{(\alpha+\beta)/2}(\x)$ допустима для $D((\alpha+\beta)/2),$ $\Gamma((\alpha+\beta)/2)),$ $Y,$ $\Delta_0$  и имеет разложение  
\begin{equation}\label{eq:CKK1}
v_{(\alpha+\beta)/2}(\x)=\sigma_q\frac{\mu_d(|\x-\y_q|)}{w_d}+\frac{c_{q}(\alpha)+c_{q}(\beta)}{2}+o(1), \, \x\to\y_q. 
\end{equation} Применяя аналог формулы \eqref{eq:ass3} (cм. доказательство  \cite[Лемма 2]{DubPril5}, \cite[Лемма 2.2]{GKP}), получим 
\begin{multline}\label{eq:ass6}
0\leq I(v_{(\alpha+\beta)/2},D((\alpha+\beta)/2)_{r})-I(h_{(\alpha+\beta)/2},D((\alpha+\beta)/2)_{r})-\\\sum_{q=1}^{l}\sigma_q\left(\frac{c_{q}(\alpha)+c_{q}(\beta)}{2}-c_q\left(\frac{\alpha+\beta}{2}\right)\right)+o(1),\,\,r\rightarrow0.
\end{multline}
Из определения функции $v_{(\alpha+\beta)/2}(\x)$ и свойства модуля вектора $|\x+\y|^2\leq 2(|\x|^2+|\y|^2)$ вытекает \begin{multline}\label{jj}
I(v_{(\alpha+\beta)/2},D((\alpha+\beta)/2)_{r})\leq \frac{1}{2}\int\limits_{D(\alpha+\beta)/2}(|\nabla( h_\alpha(\x)-h_\beta(\x^*))|^2)d\x+\frac{1}{2}\int\limits_{D(\alpha+\beta)/2}|\nabla h_\beta(\x)|^2 d\x=\\\frac{1}{2}\int\limits_{D(\alpha)}|\nabla h_\alpha(\x)|^2 d\x+\frac{1}{2}\int\limits_{D(\beta)}|\nabla h_\beta(\x)|^2 d\x.
\end{multline} Из \eqref{ass11}, \eqref{eq:ass6}, \eqref{jj} следует 
$$\sum_{q=1}^l\sigma_q c_q(\alpha)+\sum_{q=1}^l\sigma_q c_q(\beta)\leq 2\sum_{q=1}^l\sigma_q c_q\left(\frac{\alpha+\beta}{2}\right),$$
или $$ \frac{f(\alpha)+f(\beta)}{2}\leq f\left(\frac{\alpha+\beta}{2}\right).$$
Последнее неравенство и означает вогнутость  функции $f(\alpha).$  Лемма доказана.

Перейдем теперь к доказательству Теоремы 2. Заметим, что условия теоремы гарантируют, что  $\sum_{k=1}^n\delta_k=0.$
Будем считать, что $\theta_0=0$ и $\theta_m=2\pi.$ Обозначим $$B_j=D\cap \{\theta_j\leq\theta\leq\theta_{j+1}\},$$ $$ B_j^+=D\cap \{\theta_j\leq\theta\leq\frac{\theta_j+\theta_{j+1}}{2}\},\  B_j^-=D\cap \{\frac{\theta_j+\theta_{j+1}}{2}\leq \theta\leq\theta_{j+1}\},$$ $$\alpha_{j}=\frac{\theta_{j+1}-\theta_{j}}{2},$$
$j=0,\ldots,m-1.$ Пусть $Y=\{y_q\}_{q=1}^l,$ $\Delta_0=\{\sigma_q\}_{q=1}^l$ --- это точки из $X$ и соответствующие им заряды ($\sigma_q=\delta_k$ если $\y_q=\x_k$), лежащие на $\{\theta=0\}.$  Функцию $h_\alpha(\x)$ из Леммы 1, определяемую   множеством  $Y$, зарядом $\Delta_0$ и областью  $D(\alpha)=D\cap\{0<\theta<\alpha\}$, обозначим   $h^1_{\alpha}(\x).$   Аналогично  пусть $h^2_{\alpha}(\x)$  определяется набором $Y$,  $-\Delta_0=\{-\sigma_q\}_{q=1}^l$ и $D(\alpha)=D\cap\{0<\theta<\alpha\},$ $h^3_{\alpha}(\x)$ --- набором $Y$,  $-\Delta_0$ и $D(\alpha)=D\cap\{-\alpha<\theta<0\},$ и  $h^4_{\alpha}(\x)$  --- набором  $Y$,  $\Delta_0$ и $D(\alpha)=D\cap\{-\alpha<\theta<0\}.$ Константу из разложения \eqref{eq:CKK} функции $h^p_{\alpha}(\x)$  обозначим $c_q^p(\alpha),$ $p=1,2,3,4.$
 Определим функции  $$\psi_j^+(\x)=h^1_{\alpha_j}(\theta_{j}(\x)), \x\in B_j^+, j=0,2,\ldots,m-2,$$   
$$\psi_j^+(\x)=h^2_{\alpha_j}(\theta_{j}(\x)), \x\in B_j^+, j=1,3,\ldots,m-1,$$ 
$$\psi_j^-(\x)=h^3_{\alpha_j}(\theta_{j+1}(\x)), \x\in B_j^-, j=0,2,\ldots,m-2,$$   
$$\psi_j^-(\x)=h^4_{\alpha_j}(\theta_{j+1}(\x)), \x\in B_j^+, j=1,3,\ldots,m-1,$$  где обозначение $\varphi(\x)$ означает поворот на угол $\varphi$ (а именно   $\varphi(\x)=(r,\theta-\varphi,\x')$, если  $\x=(r,\theta,\x')$). В области $B_j$, $j=0,\ldots,m-1$ зададим функции 
$$\psi_j(\x)=\left\{\begin{array}{l}\psi_j^+(\x), \x\in B_j^+,\\
\psi_j^-(\x), \x\in B_j^-,\\
0,\ \x=(r, (\theta_j+\theta_{j+1})/2,x').
\end{array}\right. $$  По построению функция $\psi_j(\x)$ гармоническая в $B_j,$ имеет нулевую производную    по нормали  на границе $\partial B_j$ (за исключением точек $X$), и разложение типа \eqref{eq:BK} в окрестности точек $X\cap\overline{B_j}.$ Пусть $u(\x)$ потенциальная функция Неймана набора $X$, $\Delta$, и $\sum ^{j} \delta_{k}a_{k}$ означает суммирование тех слагаемых $\delta_k a_k$, которые соответствуют  точкам $x_k\in \overline B_j.$ Повторяя доказательство Леммы 2.2  \cite{GKP}, получим  
\begin{multline}\label{prog1}
    0\leq I(u,(B_j)_r)-I(\psi_j,(B_j)_r)-\sum ^{j} \delta_{k}a_{k}+\sum_{q=1}^l\sigma_q c_q^1(\alpha_j)+\\\sum_{q=1}^l(-\sigma_q) c_q^3(\alpha_j)+o(1),\,\, j=0,\ldots,m-2,
\end{multline}
\begin{multline}\label{prog2}
    0\leq I(u,(B_j)_r)-I(\psi_j,(B_j)_r)-\sum ^{j} \delta_{k}a_{k}+\sum_{q=1}^l(-\sigma_q) c_q^2(\alpha_j)+\sum_{q=1}^l\sigma_q c_q^4(\alpha_j)+o(1),\,\,\\  j=1,\ldots,m-1.
\end{multline}
Чтобы получить неравенство \begin{multline}\label{prog6}
   \sum _{k=1}^{n} \delta_{k}a_{k}\leq \frac{1}{2}\sum_{j=0,\ldots,m-2}\sum_{q=1}^l\sigma_q c_q^1(\alpha_j)+\frac{1}{2}\sum_{j=0,\ldots,m-2}\sum_{q=1}^l(-\sigma_q) c_q^3(\alpha_j)\\+\frac{1}{2}\sum_{j=1,\ldots,m-1}\sum_{q=1}^l(-\sigma_q) c_q^2(\alpha_j)+\frac{1}{2}\sum_{{j=1,\ldots,m-1}}\sum_{q=1}^l\sigma_q c_q^4(\alpha_j),
\end{multline} мы  просуммируем неравенства  \eqref{prog1}, \eqref{prog2} по всем $j=0,\ldots, m-1$,  применим разложение \eqref{ass2} и равенства  
  \begin{multline*}
I(\psi_j,(B_j)_r)=\sum_{q=1}^{l}\sigma_{q}^2\frac{\mu_d(r)}{w_d}+\frac{1}{2}\sum_{q=1}^l\sigma_q c_q^1(\alpha_j)+\frac{1}{2}\sum_{q=1}^l(-\sigma_q) c_q^3(\alpha_j)+o(1),\,\, j=0,\ldots,m-2,
\end{multline*}
\begin{multline*}
    I(\psi_j,(B_j)_r)=\sum_{q=1}^{l}\sigma_{q}^2\frac{\mu_d(r)}{w_d}+\frac{1}{2}\sum_{q=1}^l(-\sigma_q) c_q^2(\alpha_j)+\frac{1}{2}\sum_{q=1}^l\sigma_q c_q^4(\alpha_j)+o(1),\,\,\\  j=1,\ldots,m-1,
\end{multline*}  
$$m\sum_{q=1}^{l}\sigma_{q}^2=\sum_{k=1}^{n}\delta_k^2,$$  а также учтем   тот факт, что каждая  точка $\x_k\in X$ принадлежит двум замкнутым областям  $\overline B_j,$  Далее отметим, что из данного в Лемме 1 определения $h_\alpha(\x)$  вытекают равенства $h^1_\alpha(\x)=-h^2_\alpha(\x),$ $h^3_\alpha(\x)=-h^4_\alpha(\x),$ Следовательно, 
$$\sum_{q=1}^l\sigma_q  c_q^1(\alpha_j)=\sum_{q=1}^l(-\sigma_q) c_q^2(\alpha_j),\  \sum_{q=1}^l\sigma_q  c_q^4(\alpha_j)=\sum_{q=1}^l(-\sigma_q) c_q^3(\alpha_j).$$ Кроме того, в области $B(\alpha)=D\cap\{-\alpha<\theta<\alpha\}$ существует единственная гармоническая (за исключением точек $Y$) функция c разложением \eqref{eq:CKK} в окрестности $\y_q$, $q=1,\ldots,l,$ равная нулю на ${\partial B(\alpha)}\cap(\{\theta=\alpha\}\cup\{\theta=-\alpha\})$ и имеющая нулевую нормальную производную на оставщейся части границы ${\partial B(\alpha)}$. Указанная функция совпадает с  $h^1_\alpha(\x)$ в области $D\cap\{0<\theta<\alpha\},$  и с функцией $h^4_\alpha(\x)$ в области $D\cap\{-\alpha<\theta<0\}.$ Поэтому $c_q^1(\alpha)=c_q^4(\alpha)$ и неравенство  \eqref{prog6} принимает вид \begin{equation}\label{prog7}
   \sum _{k=1}^{n} \delta_{k}a_{k}\leq \sum_{j=0}^{m-1}\sum_{q=1}^l\sigma_q c_q^1(\alpha_j)=\sum_{j=0}^{m-1}f(\alpha_j).
\end{equation}
Из \eqref{prog7}, установленной в   Лемме 1 вогнутости функции $f(\alpha)$  и равенства  $\sum_{j=1}^m\alpha_j=\pi,$ мы  получим неравенство 
\begin{equation}\label{prog8}
   \sum _{k=1}^{n} \delta_{k}a_{k}\leq \sum_{j=0}^{m-1}f(\alpha_j)\leq mf\left(\frac{\sum_{j=0}^{m-1}\alpha_j}{m}\right)= mf\left(\frac{\pi}{m}\right).
\end{equation}

Пусть теперь  $u^*(\x)$ потенциальная функция Неймана набора $X^*$, $\Delta$ и $a_k^*$ означают соответствующие константы из асимптотического разложения.  Повторяя вышеприведенное доказательство с заменой $X$ на $X^*$ нетрудно убедиться, что во всех неравенствах выполняется знак равенства и 
\begin{equation}\label{prog8}
   \sum _{k=1}^{n} \delta_{k}a_{k}^*= mf\left(\frac{\pi}{m}\right).
\end{equation}
Таким образом, неравенство \eqref{prog8} означает  \begin{equation}\label{prog9}
   \sum _{k=1}^{n} \delta_{k}a_{k}\leq \sum _{k=1}^{n} \delta_{k}a_{k}^*.
\end{equation}
Как было отмечено в доказательстве Теоремы 1,  \eqref{prog9} эквивалентно требуемому утверждению. 

\textbf{Acknowledgements.} Работа выполнена при финансовой поддержке РФФИ (проект  20-01-00018) и  Министерства образования и науки Российской Федерации (соглашение  № 075-02-2021-1395).

\end{document}